\documentclass[12 pt]{article}
\usepackage[english]{babel}
\usepackage[active]{srcltx}
\usepackage{amsmath}
\usepackage{amsfonts,amssymb}
\usepackage[mathcal]{eucal}

\newcommand{\rav}{\stackrel{\triangle}{=}}
\newcommand{\ravref}[1]{\stackrel{(\ref{#1})}{=}}
\newcommand{\leqref}[1]{\stackrel{(\ref{#1})}{\leq}}
\newcommand{\geqref}[1]{\stackrel{(\ref{#1})}{\geq}}
\newcommand{\rref}[1]{$(\ref{#1})$}
\newcommand{\mm}[1]{{\mathbb{#1}}}
\newcommand{\av}{v}
\newcommand{\bw}{w}
\newcommand{\fl}{{*}}
\newcommand{\sr}{{*}}
\newcommand{\nr}{{*}}
\newcommand{\avm}{\mathcal{V}^{\,\fl}}
\newcommand{\avp}{\mathcal{V}^{\,\sr}}
\newcommand{\bwm}{\mathcal{W}^{\,\fl}}
\newcommand{\bwp}{\mathcal{W}^{\,\sr}}

\newcommand{\ncym}{{\widetilde{\mathbb{V}}}^{\,\!\fl}}
\newcommand{\cym}{{\mathbb{V}}^{\,\!\fl}}
\newcommand{\cyp}{{\mathbb{V}}^{\,\!\sr}}
\newcommand{\cyg}{{\mathbb{V}}}
\newcommand{\ssm}[2]{\left[#1\right]^{\,\!\fl}_{#2}}
\newcommand{\ssmk}[2]{\Big[#1\Big]^{\,\!\fl}_{#2}}
\newcommand{\ssml}[1]{\Big[#1}
\newcommand{\ssmll}[1]{\bigg[#1}
\newcommand{\ssmlll}[1]{\Bigg[#1}
\newcommand{\ssmr}[2]{#1\Big]^{\,\!\fl}_{#2}}
\newcommand{\ssmrr}[2]{#1\bigg]^{\,\!\fl}_{#2}}
\newcommand{\ssmrrr}[2]{#1\Bigg]^{\,\!\fl}_{#2}}
\newcommand{\ssp}[2]{\left[#1\right]^{\,\!\sr}_{#2}}
\newcommand{\sspk}[2]{\Big[#1\Big]^{\,\!\sr}_{#2}}

\newcommand{\sspll}[1]{\bigg[#1}
\newcommand{\ssplll}[1]{\Bigg[#1}
\newcommand{\sspr}[2]{#1\Big]^{\,\!\sr}_{#2}}
\newcommand{\ssprr}[2]{#1\bigg]^{\,\!\sr}_{#2}}
\newcommand{\ssprrr}[2]{#1\Bigg]^{\,\!\sr}_{#2}}
\newcommand{\avg}{\mathcal{V}^{\,\nr}}
\newcommand{\bwg}{\mathcal{W}^{\,\nr}}
\newcommand{\mA}{{\mathfrak{A}}}

\newcommand{\ct}[1]{{\mathcal{#1}}}
\newcommand{\epsi}{\varepsilon}
\newcommand{\bo}{\hfill {$\Box$}}

\begin{document}

\title{On an example for the Uniform Tauberian theorem in abstract control systems\thanks{
               Krasovskii Institute, Yekaterinburg, Russia;\ \  Ural Federal University, Yekaterinburg, Russia
}
}

\author{Dmitry Khlopin\\
{\it khlopin@imm.uran.ru} 
}

\maketitle

\begin{abstract}
 The paper is devoted to the asymptotic behavior of value functions of abstract control problem with
 the long-time and  discounted averages. The Uniform Tauberian Theorem for these problems states that the uniform convergence of  value functions  for long-time~averages (as the horizon tends to infinity)
  is equivalent to the uniform convergence of  value functions for discounted averages
   (as the discount tends to zero), and that the limits are identical.
   According to Miquel Oliu-Barton and Guillaume Vigeral, this assertion holds if the set of all feasible processes is closed with respect to concatenation. In this paper, we refine this condition.



\end{abstract}

 Assume the following items are given:
 \begin{itemize}
   \item   a nonempty set  $\Omega$;
   \item   a nonempty subset $\mm{K}$ of mappings from $\mm{R}_{\geq 0}$ to $\Omega$;
   \item   a running cost $g:\Omega\mapsto [0,1];$ for each process $z\in \mm{K},$ assume the  map $t\mapsto g(z(t))$ is  Borel-measurable.
 \end{itemize}
 For all $\omega\in\Omega$, define  $$\Gamma(\omega)\rav\{z\in\mm{K}\,|\,z(0)=\omega\}\quad\forall \omega\in\Omega$$
  as the set of all  feasible processes $z\in\mm{K}$ that begin at $\omega$.

 Let us now define the time average $\av_T(z)$ and the discount average $\bw_\lambda(z)$ for each process $z\in \mm{K}$
 by the rules:
 \begin{eqnarray*}
 \av_T(z)\rav\frac{1}{T}\int_{0}^T g(z(t))\,dt,\quad
 \bw_\lambda(z)\rav\lambda\int_{0}^\infty e^{-\lambda t} g(z(t))\,dt\quad \forall T,\lambda>0,z\in\mm{K}.
 \end{eqnarray*}
Also define the corresponding value functions:
\begin{eqnarray*}
 V[\av_T](\omega)\rav\inf_{z\in\Gamma(\omega)} \av_T(z),\quad
 V[\bw_\lambda](\omega)\rav\inf_{z\in\Gamma(\omega)} \bw_\lambda(z)\qquad \forall \omega\in\Omega,T,\lambda>0.
\end{eqnarray*}
   Note that the definitions are valid.

  For all processes  $z,z':\mm{R}_{\geq 0}\to\Omega$ and $h>0$ such that $z(h)=z'(0)$, define the concatenation
  $z\diamond_h z':\mm{R}_{\geq 0}\to\Omega$ is as follows:
 $(z\diamond_h z')(t)\rav
 z(t)$ if $t<h,$ and
  $(z\diamond_h z')(t)\rav z'(t-h)$,
 if $t\geq h.$

  In initial version of Theorem~6 in \cite{barton}, one stated that
  for a  set $\mm{K}$ closed with respect to concatenation (for all $\tau>0$),

 The Tauberian theorem for an abstract control systems \cite[Theorem~6]{barton} states that
  the following limits exist, are uniform in $\omega\in\Omega$, and coincide
    \begin{equation}\label{621u}
    \lim_{T\uparrow\infty} V[\av_T](\omega)
       =\lim_{\lambda\downarrow 0} V[\bw_\lambda](\omega)\quad \forall\omega\in\Omega
       \end{equation}
if at least one of these limits exists, and is uniform   in
 $\omega\in\Omega.$

Note that initially in \cite{barton} it was assumed that for this theorem to hold, the set $\mm{K}$ should be  closed with respect to concatenation. An example below will show that the requirement
``$\mm{K}$ closed with respect to concatenation''
 is insufficient
for it.

We agree that  the proof of  \cite[Theorem~6]{barton} is correct (see  the corrected version of the proof in
https://hal.archives-ouvertes.fr/hal-00661833v2) if we additionally assume that  any process
 obtained by cutting the beginning of a feasible process is also feasible, i.e.,
 for every process $z\in\mm{K}$
 and every time moment $\tau>0,$ the map $\mm{R}_{\geq 0}\ni t\mapsto z(t+\tau)$ is also in $\mm{K}$.
 In this case,
  \begin{equation}\label{621uu}
    \mm{K}=\{z'\diamond_\tau z''\,|\,\forall z',z''\in\mm{K},z'(\tau)=z''(0)\}\qquad \forall \tau>0.
       \end{equation}
 Now, we can also formulate the Uniform Tauberian Theorem for abstract control system as follows:

 {\bf Theorem.}
 {\it Assume that $\mm{K}$ satisfies \rref{621uu}.

 Then,
  all limits in \rref{621u} exist, are uniform in $\omega\in\Omega$, and coincide
if at least one of these limits exists, and is uniform   in
 $\omega\in\Omega.$
}

{\bf The example.}

Set $\Omega\rav \mm{R}_{\geq 0}\times \mm{R}_{\geq 0}\times \mm{R}_{\geq 0}. $
Define
$$g(x,y,r)\rav 0 \ \textrm{if}\ x\in[1,2],r=0\ \textrm{and}\ g(x,y,r)\rav 1,\ \textrm{otherwise}.$$
For all $\omega=(x,y,r)\in \Omega$, define $a_{\omega}:\mm{R}_{\geq 0}\to \Omega$
as follows: $$a_{\omega}(t)=a_{(x,y,r)}(t)\rav(x,y,t+r)\quad \forall t\geq 0.$$
In particular, $a_{\omega}(0)=\omega$ for all $\omega=(x,y,r)\in \Omega.$

Also, for all $s\in\mm{R}_{\geq 0}$, define $b_{s}:\mm{R}_{\geq 0}\to \Omega$
as follows:
$$b_{s}(t)\rav(st,t,0)\quad \forall t\geq 0.$$
In particular, $b_{s}(0)=(0,0,0)$ for all $s\in\mm{R}_{\geq 0}$.


Set
\begin{eqnarray*}
\mm{K}&\rav& \{a_{\omega}\,|\,\omega\in\Omega\}\cup\{b_s\,|\,s\in\mm{R}_{\geq 0}\}\cup
\{b_s\diamond_\tau a_{b_s(\tau)}\,|\,s,\tau\in\mm{R}_{\geq 0}\},\\
    \Gamma(\omega)&\rav&\{z\in\mm{K}\,|\,z(0)=\omega\}\quad\forall \omega\in\Omega.
\end{eqnarray*}
It is easy to see that $\Gamma(\omega)=\{a_\omega\}$ for all $\omega\in\Omega\setminus\{(0,0,0)\}.$
Note that  $a_{(x,y,r)}\diamond_\tau z$ is well-defined for some $z\in\mm{K}$ iff
$\tau>0$ and $z\in\Gamma(a_{(x,y,r)}(\tau))=\Gamma(x,y,r+\tau)$ iff $z=a_{a_{(x,y,r)}(\tau)}.$
Moreover, $a_{(x,y,r)}\diamond_\tau a_{a_{(x,y,r)}(\tau)}=a_{(x,y,r)}$ for all $\tau>0$.

Thus, $\mm{K}$ is closed with respect to concatenation.

It  is easy to see that $\av_T(a_\omega)=\bw_\lambda(a_\omega)=1$ for all $\lambda,T>0,\omega\in\Omega.$
Therefore,  $$V[\av_T](\omega)= V[\bw_\lambda](\omega)=1\qquad  \forall \omega\in\Omega\setminus\{(0,0,0)\}.$$

Since
$$ g(b_s(t))\leq g((b_s\diamond_\tau a_{b_s(\tau)})(t))\leq g(b_0(t))=g((b_0\diamond_\tau a_{b_0(\tau)})(t))=1 \quad \forall s,t,\tau\geq 0,$$
we obtain $$\av_T(b_s)\leq\av_T(b_s\diamond_\tau a_{b_s(\tau)}),\quad \bw_\lambda(b_s)\leq\bw_\lambda(b_s\diamond_\tau a_{b_s(\tau)}).$$
Then,
 $$V[\bw_\lambda](0,0,0)=\inf_{s> 0} \bw_\lambda(b_s),\quad
   V[\av_T](0,0,0)=\inf_{s> 0} \av_T(b_s)\qquad\forall T,\lambda>0.$$
It  is easy to prove that $\av_T(b_{T/2})=1/2\geq \av_T(b_s)$ for all $s>0$. Therefore, $$V[\av_T](0,0,0)=1/2\quad \forall T>0.$$

In addition, for any $s,\lambda>0,$ we obtain
\begin{eqnarray*}
 \bw_\lambda(b_s)&=&\lambda\int_{0}^\infty e^{-\lambda t} g(b_s(t))\,dt\\
 &=&\lambda\int_{0}^\infty e^{-\lambda t} (1-1_{[1,2]}(st))\,dt\\
 &=&1-\lambda\int_{0}^\infty e^{-\lambda t} 1_{[1,2]}(st)\,dt\\
 &=&1-\int_{0}^\infty e^{-\tau} 1_{[1,2]}(s\tau/\lambda)\,d\tau\\
 &=&1-\int_{\lambda/s}^{2\lambda/s} e^{-\tau}\,d\tau\\
 &=&1-(e^{-\lambda/s}-e^{-2\lambda/s})\geq \inf_{x\in\mm{R}}(1-x+x^2)=3/4.
\end{eqnarray*}
    On the other hand,
    for $s=\lambda/\ln 2,$ we have $e^{-\lambda/s}=1/2$, $\bw_\lambda(b_s)=1-1/2+1/4=3/4$. Hence,
    $$V[\bw_\lambda](0,0,0)=\bw_\lambda(b_{\lambda/\ln 2})=3/4\quad \forall \lambda>0.$$

    Thus,
     the following limits exist, are uniform on $\Omega$,
    \begin{equation*}
    \lim_{T\uparrow\infty}V[\av_T](\omega)=V[\av_1](\omega),\quad
    \lim_{\lambda\downarrow 0}V[\bw_\lambda](\omega)=V[\bw_1](\omega)\qquad \forall \omega\in\Omega.
       \end{equation*}

  Remember that $\mm{K}$ is closed with respect to concatenation. Then,\
  the Uniform Tauberian Theorem for abstract control system would implies that $$V[\bw_1]\equiv V[\av_1].$$
  However,
   $$3/4=V[\bw_1](0,0,0)\neq V[\av_1](0,0,0)=1/2.$$

    Therefore the Uniform Tauberian Theorem \cite[Theorem~6]{barton} does not holds for this abstract control system.


\begin{thebibliography}{99}
\bibitem{barton}
 M.~Oliu-Barton, G.~Vigeral.
\newblock  A uniform Tauberian theorem in optimal control.
In:
\newblock {\em Advances in Dynamic Games},
Birkh\"{a}user, Boston, pp. 199-215, 2013.\\
Erratum  https://hal.archives-ouvertes.fr/hal-00661833v2
\end{thebibliography}
\end{document}